\documentclass[12pt,reqno, oneside]{amsart}
\usepackage{pdfsync}  
\usepackage[width=6in, height=8.6in]{geometry}
\usepackage{hyperref} %should be loaded last

\newcommand{\seqnum}[1]{\href{http://oeis.org/#1}{{#1}}}
\newcommand{\sumz}[1]{\sum_{#1=0}^\infty}
\newcommand{\dpow}[2]{\frac{#1^{#2}}{#2!}}

\begin{document}
\title{A short proof of the Almkvist--Meurman theorem}
\author{Ira M. Gessel$^*$}
\address{Department of Mathematics\\
   Brandeis University\\
   Waltham, MA 02453}
\email{gessel@brandeis.edu}
\date{October 23, 2023}
\thanks{Supported by a grant from the Simons Foundation (\#427060, Ira Gessel).}
\begin{abstract}
We give a short generating function proof of the Almkvist--Meurman theorem:  for integers $h$ and $k\ne0$, the numbers $M_n(h,k)$ defined by 
\begin{equation*}
kx\frac{e^{hx}-1}{e^{kx}-1}=\sumz n M_n(h,k)\dpow xn
\end{equation*}
are integers. The proof is related to Postnikov's functional equation for the generating function for intransitive trees.
\end{abstract}

\maketitle
\thispagestyle{empty}
\subsection{Introduction}
The Bernoulli polynomials $B_n(u)$ may be defined by 
\begin{equation}
\label{e-bp}
\sum_{n=0}^\infty B_n(u) \frac{x^n}{n!}=\frac{xe^{ux}}{e^x-1},
\end{equation}
and the $n$th Bernoulli number $B_n$ is $B_n(0)$. 

For $k\ne0$, let
\begin{equation}
\label{e-Mdef}
M_n(h,k)=k^n\bigl(B_n(h/k) - B_n\bigr), 
\end{equation}
so 
\begin{equation}
\label{e-amgf}
\sumz n M_n(h,k)\dpow xn
=
kx\frac{e^{hx}-1}{e^{kx}-1}
\end{equation}

Almkvist and Meurman  \cite{am} showed  that if $h$ and $k$ are integers then $M_n(h,k)$ is an integer. Other proofs were given by Sury \cite{sury},
 Bartz and Rutkowski \cite{br},  Clarke and Slavutskii \cite{cs}, and the author \cite{gessel}. The Almkvist--Meurman  theorem also follows immediately from a much earlier result of Vandiver \cite{vandiver1937}, as explained in \cite{gessel}. Additional historical remarks can also be found in \cite{gessel}. The case $h=1$ of the theorem (from which the general case is easily derived, as explained below) was recently rediscovered by Farhi \cite{farhi22}, who gave a new proof. 
 
 We give here another  new  proof, using generating functions, of the Almkvist--Meurman theorem. The proof in  \cite{gessel} also uses generating functions, but in a different way.

A \emph{Hurwitz series} is a (formal) power series $\sum_{n=0}^\infty f_n x^n/n!$ for which each $f_n$ is an integer; i.e., it is the exponential generating function for a sequence of integers.  It is well known that Hurwitz series are closed under addition and multiplication.
Moreover, if $f$ and $g=\sum_{n=1}^\infty g_n x^n/n!$ are Hurwitz series then the composition $f(g(x))$ is a Hurwitz series, and if $g_1=1$ then $g(x)$ has a compositional inverse that is a Hurwitz series.

We will define a \emph{generalized Hurwitz series} to be a power series whose coefficients, as an exponential generating function, are polynomials in some set of variables (in our application $\alpha_1$, $\alpha_2$, 
$\beta_1$, $\beta_2$) with integer coefficients. 

Since all of the generating functions we will be concerned with are exponential, by the ``coefficients" of a generating function we mean its coefficients as an exponential generating function. 

\subsection{Proof of the Almkvist--Meurman theorem}
We first note that the general case of the Almkvist--Meurman theorem can be derived from the case in which $h=1$ and $k$ is positive, since for any integers $h$ and $k$ (with $k\ne0$), $kx (e^{hx} -1)/(e^{kx}-1)$ can be expressed as a Hurwitz series times $|k|x (e^x-1)/(e^{|k|x}-1)$.
This follows easily from 
$e^{hx}-1= -e^{hx}(e^{-hx}-1)$, $e^{kx}-1 = -e^{kx} (e^{-kx}-1)$, and for $h$ a positive integer,
$e^{hx}-1= (1+e^x +\cdots+e^{(h-1)x})(e^x-1)$.
When $h=1$ and $k$ is a positive integer, the generating function \eqref{e-amgf} may be written
\begin{equation}
\label{e-amh1}
\frac{kx}{1+e^x+\cdots +e^{(k-1)x}}.
\end{equation}

To motivate our proof, we first prove the case $h=1,k=2$ of the Almkvist--Meurman theorem, for which the generating function \eqref{e-amh1} is $2x/(e^x+1)$,
the generating function for the Genocchi numbers,  which is sequence \seqnum{A001469} in the OEIS (On-Line Encyclopedia of Integer Sequences \cite{oeis}).

We start by looking at what seems to be a completely unrelated sequence of integers, sequence \seqnum{A007889}, which has several combinatorial interpretations. These numbers count alternating trees, local binary search trees, and regions of the Linial arrangement.
Postnikov \cite{postnikov} showed that the exponential generating function $A=A(x)$ for these numbers (where we take the constant term to be 0), satisfies
\begin{equation}
\label{e-A1}
1+A = e^{(x/2)(2+A)}.
\end{equation}
In this form it is not clear that $A$ is a Hurwitz series, though it is clear that $2$ is the only prime that can occur in the denominators of the coefficients of $A$. However we can rearrange \eqref{e-A1} in a way that does show that $A$ is a Hurwitz series. 
Squaring both sides of \eqref{e-A1} gives
\begin{equation}
\label{e-1.3}
(1+A)^2=1+A(2+A) = e^{x(2+A)}.
\end{equation}
Subtracting 1 from each side and dividing by $2+A$ gives
\begin{equation}
\label{e-A1.5}
A = \sum_{n=1}^\infty (2+A)^{n-1}\dpow xn,
\end{equation}
from which it is clear that $A$ is a Hurwitz series.

Solving \eqref{e-A1} for $x$ shows that the compositional inverse of $A(x)$ is 
\begin{equation}
\label{e-A2}
\frac{2\log(1+x)}{2+x},
\end{equation}
which is therefore a Hurwitz series. Although we will not need this fact, the coefficient of $x^n/n!$ in \eqref{e-A2} is $(-1)^{n-1}\sum_{i=0}^{n-1} i!\, (n-i-1)!$.%
\footnote{This can be proved easily from the beta integral: 
%since
%\begin{equation*}
%\int_0^1 u^i(1-u)^{n-i-1}\, du = \frac{i!\,(n-i-1)!}{n!},
%\end{equation*}
%we have 
%\begin{align*}
%\sum_{n=1}^\infty (-1)^{n-1}\sum_{i=0}^{n-1} i!\, (n-i-1)! \dpow xn
%  &=\int_0^1 \sum_{n=1}^\infty (-1)^{n-1}\sum_{i=0}^{n-1} u^i (1-u)^{n-i-1}x^n\,du\\
%  &=\int_0^1 \sum_{i,j=0}^\infty (-1)^{i+j} u^i (1-u)^j x^{i+j+1}\,du\\
%  &=x\int_0^1  \frac{du}{(1+ux)(1+(1-u)x)}=\frac{2\log(1+x)}{2+x}.
%\end{align*}
%
since
\begin{equation*}
\int_0^1 z^i(1-z)^{j}\, dz = \frac{i!\,j!}{(i+j+1)!},
\end{equation*}
we have 
\begin{align*}
\sum_{i,j=0}^\infty \frac{i!\, j!}{(i+j+1)!} u^i v^j
  &=\int_0^1 \sum_{i,j=0}^\infty (uz)^i (v(1-z))^j \, dz   \\
  &=\int_0^1 \frac{dz}{(1-uz)(1-v(1-z))}\\
  &= \frac{\log(1-u) +\log(1-v)}{uv-u-v}.
\end{align*}
Setting $u=v=-x$ and multiplying by $x$ gives the formula for the coefficients of \eqref{e-A2}.
}%end of footnote
These numbers are \seqnum{A003149} in the OEIS.

Replacing $x$ with $e^x-1$ in \eqref{e-A2} gives another Hurwitz series,
\begin{equation*}
\frac{2x}{e^x+1},
\end{equation*}
the exponential generating for the Genocchi numbers, thus proving the case $k=2$ of the Almkvist--Meurman theorem.

A related proof of the integrality of the Genocchi numbers was given by Farhi \cite{farhi22a}, who used the fact that
\eqref{e-A2} is a Hurwitz series, but did not use \eqref{e-A1.5} or the connection with the alternating tree numbers.

To extend this proof to the case $h=1$, $k$ a positive integer of the Almkvist--Meurman theorem, which, as we have seen, implies the general theorem,  we replace the left side of \eqref{e-1.3} with
\begin{equation*}
(1+A)^k = 1+A\left[\binom{k}{1} + \binom{k}{2}A + \cdots + \binom{k}{k} A^{k-1}\right]
 =1+Ap_k(A),
\end{equation*}
where $p_k(u)$ is the polynomial 
\[\binom{k}{1} + \binom{k}{2}u + \cdots + \binom{k}{k} u^{k-1}=\bigl((1+u)^k -1\bigr)/u.\]
We then replace the right side of \eqref{e-1.3} with $e^{xp_k(A)}$, so we consider the functional equation
\begin{equation}
\label{e-kA}
(1+A)^k = e^{xp_k(A)}.
\end{equation}
Subtracting $1$ from each side of \eqref{e-kA} and dividing by $p_k(A)$ gives 
the functional equation
\begin{equation}
\label{e-Af}
A = \sum_{n=0}^\infty p_k(A)^{n-1}\dpow xn,
\end{equation}
which shows that $A$ is a Hurwitz series.

Solving \eqref{e-kA} for $x$ shows that the compositional inverse of $A$ is
\begin{equation}
\label{e-iAk}
\frac{k\log (1+x)}{p_k(x)}=\frac{kx\log(1+x)}{(1+x)^k-1}
\end{equation}
 so \eqref{e-iAk} is a Hurwitz series. Replacing $x$ with $e^x-1$ in \eqref{e-iAk} gives
 \begin{equation*}
kx\frac{e^x-1}{e^{kx}-1},
\end{equation*}
which must therefore be a Hurwitz series. This completes our proof of the Almkvist--Meurman theorem.

\subsection{Additional comments}
I do not know of a formula for the coefficients of \eqref{e-iAk} analogous to the formula for the coefficients of \eqref{e-A2}, nor do I know of a combinatorial interpretation to \eqref{e-kA} or \eqref{e-kA2} that generalizes the combinatorial interpretations for the case $k=2$. If we set $B=1+A$ then \eqref{e-kA} may be written in the somewhat more suggestive form
\begin{equation}
\label{e-kA2}
B = e^{x(1+B+\cdots+B^{k-1})/k}.
\end{equation}

There is an interesting generalization of the fact that $2\log(1+x)/(2+x)$ is a Hurwitz series whose compositional inverse counts certain trees, although in this generalization there is no analogue of the replacement of $x$ with $e^x-1$ in \eqref{e-A2} that gave us the the Genocchi numbers. 

Consider the functional  equation for $F$,
\begin{equation}
\label{e-drake1}
\frac{(1+\alpha_1 F)(1+\beta_2 F)}{(1+\alpha_2 F)(1+\beta_1 F)}=e^{((\alpha_1\beta_2 -\beta_1\alpha_2)F +\alpha_1-\beta_1 -\alpha_2 +\beta_2)x}.
\end{equation}
If we set $\alpha_1=\beta_2 = 1$ and $\alpha_2=\beta_1=0$ then  \eqref{e-drake1} reduces to \eqref{e-1.3}, with $A$ replacing $F$.
Solving \eqref{e-drake1} for $x$ shows that 
$F$ is the compositional inverse of 
\begin{equation}
\label{e-drake2}
\frac{1}{(\alpha_1\beta_2 -\beta_1\alpha_2)x + \alpha_1 -\beta_1 -\alpha_2+\beta_2}
  \log\left(\frac{(1+\alpha_1 x)(1+\beta_2 x)}{(1+\alpha_2 x)(1+\beta_1 x)}\right).
\end{equation}
It is not obvious from \eqref{e-drake1} or \eqref{e-drake2} that $F$ is a generalized Hurwitz series
but this follows from a result of 
Drake \cite[Proposition 1.8.1 and Theorem 1.8.4]{drake} who showed that the coefficient of $x^n/n!$ in \eqref{e-drake2} is a polynomial in $\alpha_1, \alpha_2, \beta_1, \beta_2$ homogeneous of degree $n-1$, and the coefficient of 
$\alpha_1^{a_1}\alpha_2^{a_2}\beta_1^{d_1}\beta_2^{d_2}x^n/n!$ in \eqref{e-drake2}, where $n=a_1+a_2+d_1+d_2+1$, is 
$(-1)^{n-1} (a_1+a_2)!\, (d_1+d_2)! \binom{a_1+d_1}{a_1}\binom{a_2+d_2}{a_2}$. This generalizes the formula given earlier for the coefficients of \eqref{e-A2}. That $F$ is a generalized Hurwitz series also follows from the combinatorial interpretation for its coefficients (see Drake \cite{drake}  and Gessel, Griffin, and Tewari \cite{ggt}, and the references cited there).

We can also rewrite \eqref{e-drake1} in a form analogous to \eqref{e-A1.5}.
We first subtract 1, from \eqref{e-drake1}, getting
\begin{multline*}\quad
\frac{F}{(1+\beta_1 F)(1+\alpha_2 F)}\cdot\bigl((\alpha_1\beta_2 - \alpha_2\beta_1)F + \alpha_1 -\alpha_2-\beta_1+\beta_2\bigr)\\
=\sum_{n=1}^\infty \bigl((\alpha_1\beta_2-\alpha_2\beta_1)F+\alpha_1-\alpha_2-\beta_1+\beta_2\bigr)^n \dpow xn,\quad
\end{multline*}
from which we obtain
\begin{equation}
\label{e-ss}
F = (1+\beta_1 F)(1+\alpha_2 F)\sum_{n=1}^\infty \bigl((\alpha_1\beta_2-\alpha_2\beta_1)F+\alpha_1-\alpha_2-\beta_1+\beta_2\bigr)^{n-1} \dpow xn
\end{equation}
which makes it clear that $F$ is a generalized Hurwitz series.

Unfortunately it does not seem possible to get a generalization of the Genocchi numbers from \eqref{e-drake2}.

%for BibTex
%\bibliography{bibfilename}{amtbib}

\begin{thebibliography}{99}
\bibitem{am}
G. Almkvist and A. Meurman, Values of Bernoulli polynomials and Hurwitz's zeta function 
at rational points, C. R. Math. Rep. Acad. Sci. Canada 13 (1991), 104--108.

\bibitem{br}
K. Bartz and J. Rutkowski,
On the von Staudt--Clausen theorem,
C. R. Math. Rep. Acad. Sci. Canada 15 (1993), no. 1, 46--48.

%\bibitem{carlitz}
%L. Carlitz, The Staudt--Clausen theorem, Mathematics Magazine 34 (1961), 131--146.

\bibitem{cs}
F. Clarke and I. Sh. Slavutskii,
The integrality of the values of Bernoulli polynomials and of generalised Bernoulli numbers,
Bull. London Math. Soc. 29 (1997), 22--24.

\bibitem{drake}
B. Drake, An inversion theorem for labeled trees and some limits of areas under lattice paths, Ph.D. thesis, Brandeis University, 2008.

\bibitem{farhi22a}
B. Farhi, The integrality of the Genocchi numbers obtained through a new identity and other results,
Notes on Number Theory and Discrete Mathematics 
28 (2022), no. 4, 749--757.

\bibitem{farhi22}
B. Farhi, A new generalization of the Genocchi numbers and its consequence on the Bernoulli polynomials,
Advances in Pure and Applied Mathematics 13 (2022), no. 4, 18--28.
 

\bibitem{gessel}
I. M. Gessel,
On the Almkvist--Meurman theorem for Bernoulli polynomials, 
Integers 23 (2023), Paper No. A14, 10 pp.
\bibitem{ggt}
I. M. Gessel, S. T. Griffin, and V. Tewari, Labeled binary trees, subarrangements of the Catalan arrangements, and Schur positivity, Adv. Math. 356 (2019), 106814, 67 pp.

\bibitem{oeis}
OEIS Foundation Inc., The On-Line Encyclopedia of Integer Sequences, Published electronically at \url{https://oeis.org},  2023.

\bibitem{postnikov}
A. Postnikov, Intransitive trees, J. Combin. Theory Ser. A 79 (1997), 360--366.

\bibitem{sury}
B. Sury, The value of Bernoulli polynomials at rational numbers, Bull. London Math. Soc. 25 (1993), 327--329.


\bibitem{vandiver1937}
H. S. Vandiver, 
On generalizations of the numbers of Bernoulli and Euler, 
Proc. Nat. Acad. Sci. U.S.A. 23 (1937), 555--559.

%\bibitem{vandiver1941}
%H. S. Vandiver, Simple explicit expressions for generalized Bernoulli numbers of the first order, Duke Math. J. 8 (1941), 575--584.
%

\end{thebibliography}
%\bibliographystyle{amsplain}

\end{document}